\documentclass[12pt,thmsa,emstex]{article}
\usepackage{amsfonts}
\usepackage{amssymb}
\usepackage{amsmath}
\usepackage{latexsym}
\usepackage{verbatim}
\usepackage{url}

\newtheorem{theorem}{Theorem}[section]
\newtheorem{lemma}[theorem]{Lemma}
\newtheorem{corollary}[theorem]{Corollary}
\newtheorem{proposition}[theorem]{Proposition}

\newtheorem{exercise}[theorem]{Exercise}

\newtheorem{remark}{Remark}
\newtheorem{example}{\bf{Example}}

\def\bit{\begin{itemize}}
\def\eit{\end{itemize}}

\reversemarginpar   

\def\bc{\begin{center}}
\def\ec{\end{center}}

\def\bthm{\begin{theorem}}
\def\ethm{\end{theorem}}

\def\bcor{\begin{corollary}}
\def\ecor{\end{corollary}}

\def\bprop{\begin{proposition}}
\def\eprop{\end{proposition}}

\def\blem{\begin{lemma}}
\def\elem{\end{lemma}}

\def\bex{\begin{example}}
\def\eex{\end{example} \hfill $\diamond$ \\[2ex]}

\def\bexo{\begin{exercise}}
\def\eexo{\end{exercise} }

\def\brem{\begin{remark}}
\def\erem{\end{remark}}
\def\prf{{\bf Proof: }}

\def\bdes{\begin{description}}
\def\edes{\end{description}}

\def\ita{\item[(a)]}
\def\itb{\item[(b)]}

\def\itc{\item[(c)]}

\def\iti{\item[(i)]}
\def\itii{\item[(ii)]}

\def\itiii{\item[(iii)]}
\def\itiv{\item[(iv)]}

\def\beq{\begin{equation}}
\def\eeq{\end{equation}}

\def\ben{\begin{enumerate}}
\def\een{\end{enumerate}}

\def\beqar{\begin{eqnarray}}
\def\eeqar{\end{eqnarray}}

\def\beqarr{\begin{eqnarray*}}
\def\eeqarr{\end{eqnarray*}}

\def\qed{\hfill $\Box$ \\[2ex]}
\def\prf{{\bf Proof: }\hspace{.1in}}

\catcode`\@=11

\newcommand{\F}{\mathcal{F}}

\newcommand{\PR}{\mathcal{P}}

\def\Ind{{\mathbf 1}}
\def\RR{{\mathbb R}}  

\def\NN{{\mathbb N}}

\def\Pr{{\mathsf P}}
\def\cF{{\cal F}}

\def\rar{\rightarrow}
\def\eps{\varepsilon}

\begin{document}


\title{Smale Strategies for Network Prisoner's Dilemma Games}

\author{ Kashi Behrstock\thanks{Kashi Abhyankar  Behrstock  works in the financial industry in New York, USA.}, 
Michel Bena\"{i}m\thanks{Institut de Math\'{e}matiques, Universit\'{e} de Neuch\^{a}tel,  Switzerland.} 
and Morris W Hirsch\thanks{Mathematics Department, University Wisconsin at Madison and University of California at Berkeley, USA.} }
\maketitle

\begin{abstract}
Smale's approach \cite{Smale80} to the classical two-players repeated Prisoner's  Dilemma game is revisited here for $N$-players and Network games in the framework of Blackwell's approachability, stochastic approximations and differential inclusions.
\end{abstract}
\paragraph{Keywords: Repeated prisoner dilemma, approachability,
  stochastic approximation, differential inclusions
}
{\small \tableofcontents}
\newpage
\section{Introduction}
It has been well known for many years that mutual cooperation is a
Nash equilibrium outcome in a two players infinitely  repeated
Prisoner's Dilemma game, even though defection is the dominant
strategy of the one-shot game (see e.g~ the classical book by Axelrod
\cite{AX84}).

In 1980 Smale \cite{Smale80} studied the two players repeated
Prisoner's Dilemma game under the assumption that both players have
limited memory and only keep track of the cumulative average
payoffs. In this setting, he showed that a very simple deterministic
strategy called a {\em good strategy}, if adopted by one player, leads
to cooperation, in the sense that the other player has interest to
cooperate. A good strategy, as defined by Smale, is a strategy such
that the player cooperates unless her average payoff to date is
significantly less than her opponent.  Later, Bena\"{i}m and Hirsch
\cite{BH96a} considered the stochastic analogue of Smale's
solution. In 2005, Bena\"{i}m Hofbauer and Sorin \cite{BHS05} using
tools from stochastic approximation and differential inclusions showed
that the results of Smale, Bena\"{i}m and Hirsch can be reinterpreted
in the framework of Blackwell's approachability theory
\cite{Blackwell56}, and that the assumption that "both" players keep
track only of the cumulative average payoff is unnecessary.

The present paper extends these works to variant of the classical
Prisoner Dilemma game including $N$-Players where the underlying
structure is a network. It is based on K. Abhyankar's PhD thesis
\cite{Abhyankar01}, Blackwell's approachability \cite{Blackwell56} and
the stochastic approximation approach to differential inclusions
developed in \cite{BHS05}.

Section \ref{sec:notation} sets up the notation and reviews briefly
Blackwell's approachability and some of the results in \cite{BHS05}.
Section \ref{sec:NPD} considers $N$-players prisoner dilemma games and
Section \ref{sec:graphPD} prisoner dilemma games in which players are
located at the vertices of a symmetric graph and interact only with their
neighbors. Smale good strategies are defined for these games and are
shown to be Nash equilibria.

\section{Notation and Background} \label{sec:notation}
Let $A$ and $B$ be two finite sets representing respectively the {\em
  action sets} of some {\em decision maker} $DM$ (for instance a player, or
a group of players) and the action set of {\em Nature} (for instance
the player's opponents).  Let $U: A \times B \mapsto \RR^N$ be a
vector valued {\em payoff function}.

Throughout, we let $E \subset \RR^N$ denote the convex hull of the payoff vectors
$$E = \mathsf{conv}\{U(a,b) \: : a \in A, b \in B\}.$$

At discrete times $n = 1,2, \ldots,$ $DM$ and Nature choose their
actions $(a_n,b_n) \in A \times B.$ We assume that: \bdes \ita The
sequence $\{(a_n,b_n)\}_{n \geq 0}$ is a random process defined on
some probability space $(\Omega, \cF, \Pr)$ and adapted to some
filtration $\{\cF_n\}$ (i.e~ $\{\cF_n\}$ is an increasing family of
sub-$\sigma$ fields of $\cF,$ and for each $n$ $(a_n,b_n)$ is $\cF_n$
measurable).  Here $\cF_n$ has to be understood as the history up to
time $n.$ \itb Given the history $\cF_n,$ $DM$ and Nature act
independently:
$$\Pr((a_{n+1}, b_{n+1}) = (a,b) | \cF_n) =
\Pr(a_{n+1} = a |  \cF_n) \Pr(b_{n+1} = b |  \cF_n).$$
\edes
Let $\PR(A)$ (respectively $\PR(B)$) denote the set of probabilities
over $A$ (respectively, $B$.

 A (long term) {\em strategy} for $DM$ is a stochastic process $\Theta
 = \{\Theta_n\}$ adapted to $\{\cF_n\}$ taking values in $\PR(A).$ We
 say that $DM$ uses strategy $\Theta$ if \beq
\label{eq:defstrat}
\Theta_n(a) = \Pr(a_{n+1} = a |  \cF_n)
\eeq
for all $a \in A.$

The {\em cumulative average payoff} at time $n$ is the vector
\beq
\label{eq:defun}
u_n = \frac{1}{n} \sum_{k = 1}^n U(a_k,b_k) \in E.
\eeq
 Strategy $\Theta$ is said to be  {\em payoff-based} provided
$$\Theta_n(a)  = Q_{u_n}(a)$$
 for all $a \in A,$ where for each $u \in E$, $Q_u( \cdot)$ is a
 probability over $A$ and $u \in E \mapsto Q_u \in \PR(A)$ is
 measurable.  In this case, the family $Q = \{Q_u\}_{u \in E}$ is
 identified with $DM$'s {\em strategy}.

\bex[$M$-Players games] \label{ex:Nplayers}
{\rm Consider an $M$-players game with $M \geq 2.$  Players are denoted
  $i = 1, \ldots, M.$  Player $i$ has a finite {\em action set} (or
  {\em pure
  strategy set}) denoted $\Sigma^i$, and a
{\em payoff function}
$U^i : \Sigma^1 \times \ldots \times \Sigma^M \mapsto \RR^{N_i}$ for
some $N_i \geq 1.$

At each discrete time $n = 1, 2, \ldots$ Player $i$ chooses an action
$s^i_n \in \Sigma^i$ and receives the payoff $U^i(s^1_n, \ldots,
s^M_n).$

Choose $DM$ to be some given subset of players, say $I = \{1 , \ldots,
k\}.$ Set $$A = \Sigma^1 \times \ldots \times \Sigma^k, B =
\Sigma^{k+1} \times \ldots \times \Sigma^M$$ and $U : \Sigma \mapsto \RR^{N_1} \times \ldots \times \RR^{N_2} \simeq
\RR^N, N = \sum_{i = 1}^M N_i,$ with
 $$U(a,b) = (U^1(a,b), \ldots U^M(a,b)).$$
}
\eex
\subsection*{The Limit Set theorem}
Assume that $DM$ has a payoff-based strategy $Q$. For each $u \in E$
let
\beq
\label{eq:defC}
C(u) = \left \{ \sum_{a \in A, b \in B} U(a,b) Q_u(a) \nu(b) \: : \nu
\in \PR(B) \right \}.  \eeq
The set $C(u)$ is the convex set
containing all the average payoffs that are obtained when $DM$ plays
the mixed strategy $Q_u$ and Nature plays any mixed strategy.

Let $\mathbf{C} \subset E \times E$ be the intersection of all closed
subset $\mathbf{G} \subset E \times E$ for which the fiber
 $\{y \in E \: : (x,y) \in \mathbf{G}\}$ is convex and contains $C(x)$. The {\em
  closed-convex} extension of $C,$ denoted $\overline{\mathsf{co}}(C)$
is defined as
$$\overline{\mathsf{co}}(C)(x) = \{y \in E \: : (x,y)
\in\mathbf{C}\}.$$ For convenience we extend
$\overline{\mathsf{co}}(C)$ to a set-valued map
$\overline{\mathsf{co}}(C)$ on $\RR^N$, also denoted
$\overline{\mathsf{co}}(C)$, by setting
\beq
\label{eq:defChat}
\overline{\mathsf{co}}(C)(x) =\overline{\mathsf{co}}(C)(r(x)).
\eeq
where for all $x \in \RR^N, r(x) \in E$ denotes the unique point in $E$ closest to $x.$
Associated to $\overline{\mathsf{co}}(C)$  is
 the differential inclusion
\beq
\label{eq:inclusi}
\frac{du}{dt} \in F(u) := - u + \overline{\mathsf{co}}(C)(u).
\eeq
A {\em solution} to (\ref{eq:inclusi}) is an absolutely continuous
mapping $t \mapsto \eta(t)$ verifying $\dot{\eta}(t)  \in F(\eta(t))$ for
almost every $t \in \RR.$ Given such a solution, its {\em initial
  condition} is the point $\eta(0).$
Throughout, we let $S_u  \subset C^0(\RR, \RR^N)$ denote the set of
all solutions to (\ref{eq:inclusi}) with initial condition $u.$
By construction, $F$ maps points to non empty compact convex sets and
has a closed graph. Thus, by standard results on differential
inclusions, $S_u$ is a nonempty subset of $C^0(\RR, \RR^N)$ htat
iscompact (for the topology of uniform convergence on compact
intervals) and (\ref{eq:inclusi}) induces a set-valued dynamical
system $\Phi = \{\Phi_t\}$ defined for all $t \in \RR$ and $u \in
\RR^N$ by
$$\Phi_t(u) = \{ \eta(t) \: : \eta \in S_u \}.$$
A set $\Lambda \subset \RR^N$ is said to be {\em invariant} for
(\ref{eq:inclusi}) if for all $u \in \Lambda$ there exists $\eta \in
S_u$ such that $\eta(\RR) \subset A$ (see section 3 of \cite{BHS05}
for other notions of invariance, more details and references on set
valued dynamics).

A nonempty compact set $\Lambda$ is called an {\em attracting set} for
$\Phi$ provided there is some neighborhood $U$ of $\Lambda,$ called a
{\em fundamental neighborhood}, with the property that for every $\eps
> 0$ there exists $t_\eps > 0$ such that $\Phi_t(U) \subset
N^{\eps}(\Lambda)$ for all $t \geq t_{\eps}.$ Here $N^{\eps}$ stands
for the $\eps$ neighborhood of $\Lambda.$ If in addition $\Lambda$ is
invariant, $\Lambda$ is called an {\em attractor}.  By Proposition
3.10 in \cite{BHS05}, every attracting set contains an attractor with
the same fundamental neighborhood.

The {\em basin of attraction} of an attracting set $\Lambda$ is the set
$${\cal W}(\Lambda) = \{u \in \RR^N \: : \omega_{\Phi}(u) \subset \Lambda\}$$
where $$\omega_{\Phi}(u) = \bigcap_{t \geq 0} \overline{\Phi_{[t,\infty[}(u)}.$$
We let
$$L = L(\{u_n\})$$ denote the limit set of the sequence $\{u_n\}$
defined by (\ref{eq:defun}). Note that $L$ is a random subset of $E.$

Point $p \in \RR^N$ is called {\em attainable} if for any $n \in \NN$
and any neighborhood $U$ of $p$
$$\Pr(\exists m \geq n :\: u_m \in U) > 0.$$ We let
$\mathrm{Att}(\{u_n\})$ denote the set of attainable points.

Parts $(i)$ and $(ii)$ of the following result follow from Theorems
3.6 and 3.23 in \cite{BHS05}, generalizing the limit set theorem
obtained for stochastic approximation processes (associated to an ODE)
in \cite{B96, B99} and asymptotic pseudotrajectories (of an ODE) in
\cite{BH96}. Part $(iii)$ follows from \cite{FR10} generalizing a
result obtained for stochastic approximation processes (associated to
an ODE) in \cite{B99}.
  \bthm
\label{th:limitset}
Suppose that $DM$  uses the payoff-based strategy $Q.$ Then with
probability one (regardless of Nature  strategy)
\bdes

\iti  $L = L(\{u_n\})$ is almost surely an internally chain-transitive
set of (\ref{eq:inclusi}).

\itii If $\Lambda$ is an attracting set for $\Phi$ then $L \subset
\Lambda$ on the event $L \cap {\cal W}(\Lambda) \neq \emptyset,$

\itiii If $\mathrm{Att}(\{u_n\}) \cap {\cal W}(\Lambda) \neq
\emptyset$ then $\Pr(L \subset \Lambda) > 0.$

\edes
\ethm
We refer the reader to \cite{BHS05} for the definition of "internally
chain-transitive" sets since this notion will not be used here but for
the fact that an internally chain-transitive set is compact and
invariant under differential inclusions (\ref{eq:inclusi}).

\subsection*{Approachability}
Let $d$ denote the Euclidean distance in $\RR^N.$ A set $\Lambda \subset E$ is said {\em approachable} if there exists a long term strategy for $DM$ such that, regardless of Nature strategy, $$d(u_n, \Lambda) \rar 0.$$
Given a compact
subset $\Lambda \subset E$ and $x \in E$, define
\[
\Pi_{\Lambda}(x) = \{y \in \Lambda : d(x,\Lambda) = d(x,y)\}
\]
where $d(x,\Lambda) = \inf \{d(x,y) \: : y \in \Lambda\}.$

 Record that $N^r(\Lambda) = \{x \in E \: : d(x,\Lambda) < r\}.$
We say that $\Lambda$ is a {\em local } {\em ${\cal B}$-set} for the
payoff-based strategy $Q$ (or simply a  local ${\cal B}$-set) if there
exists $r > 0$ such that for all $x \in N^r(\Lambda) \setminus
\Lambda$ there exists $y \in \Pi_{\Lambda}(x)$ such that the
hyperplane orthogonal to $[x,y]$ at $y$ separates $x$ from
$C(x)$. That is,
\beq
\label{eq:separ}
\langle x-y, v-y\rangle \leq 0
\eeq
 for all $v \in C(x)$ as defined by (\ref{eq:defC}).
If $\Lambda$ is a local ${\cal B}$-set for all $r > 0$ it is  simply called a {\em ${\cal B}$-set.}
Blackwell \cite{Blackwell56}, proved that being a ${\cal B}$-set is a
sufficient condition for approachability.

\bthm
\label{th:Black}
Let $\Lambda \subset E$ be a local ${\cal B}$-set for the payoff-based
strategy $Q.$ Then
\bdes

\iti $\Lambda$ contains an attractor for
$\Phi$ with fundamental neighborhood $U = N^r(\Lambda).$ In
particular,
\bdes \ita $L \subset \Lambda$ on the event $L \cap U \neq
\emptyset.$ \itb $\mathrm{Att}(\{u_n\}) \cap U \neq \emptyset
\Rightarrow \Pr(\{L \subset \Lambda\}) > 0.$ \edes \itii If $\Lambda$
is a B-set, then $\Pr(L \subset \Lambda) = 1.$ \edes \ethm \prf It is
proved in \cite{BHS05}, Corollary 5.1 that $\Lambda$ contains an
attractor for (\ref{eq:inclusi}) provided inequality (\ref{eq:separ})
holds for all $v \in \overline{\mathsf{co}}(C)(x)$ (rather than merely
$v \in C(x)$).  It then suffices to prove that (\ref{eq:separ}) also
holds for all $v \in \overline{\mathsf{co}}(C)(x).$

Let
$\mathrm{Graph}(C) = \{(x,y) \in E \times E \:: y \in C(x)\}$.  Denote
its closure  by
$\overline{\mathrm{Graph}}(C)$, and set
\[\overline{\mathrm{Graph}}_x(C) = \{y \in E \:: (x,y) \in
\overline{\mathrm{Graph}}(C)\}.
\]
Let $D(x)$ be the convex hull of
 $\overline{\mathrm{Graph}}_x(C).$ It follows from (\ref{eq:separ})
 and compactness of $\Lambda$ that $\langle x - y, v - y \rangle \leq
 0$ for all $x \in E, v \in \overline{\mathrm{Graph}}_x(C)$ and some
 $y \in \Pi_{\Lambda}(x).$ Clearly, this inequality still holds for
 all $v \in D(x).$ We claim that $\overline{\mathsf{co}}(C)(x) = D(x)$
 from which the proof of $(i)$ follows.

  Proof of the claim: The inclusion $D(x) \subset
  \overline{\mathsf{co}}(C)(x)$ follows from the definitions. To prove
  the opposite inclusion it suffices to verify that
 $\mathrm{Graph}(D)$ is
  closed. Let $x_n \rar x$, \,$ y_n  \rar y$ with $y_n \in D(x_n).$

By the  Caratheodory Theorem (see e.g~ Theorem 11.1.8.6 in \cite{Ber78}),
 the convex hull of a set $G \subset \RR^N$ equates the set obtained by taking all convex combinations of $N+1$ points in $G.$ Thus, there exist
\[w_n = (w_{n,1}, \ldots, w_{n,N+1}) \in
\overline{\mathrm{Graph}}_{x_n}(C)^{N+1}
\]
 and
\[\alpha_n =
(\alpha_{n,1}, \ldots, \alpha_{n,N+1}) \in \Delta^N \ \text {(the unit
$N$-dimensional simplex of $\RR^{N+1}$)}
\] such that $$y_n = \sum_{i =   1}^{N+1} \alpha_{n,i} w_{n,i}.$$

By
 compactness, after replacing sequences by subsequences we can assume
 that $\alpha_n \rar \alpha \in \Delta^N$ and $w_n \rar w.$
 Closedness of  $ \overline{\mathrm{Graph}}(C)$ ensures that $w \in
 \overline{\mathrm{Graph}}_{x}(C)^{N+1}.$ Thus $y \in D(x).$ This
 proves the claim.

Assertions $(a)$ and $(b)$ are now  consequences of Theorem
\ref{th:limitset}. The lastr  statement was proved by Blackwell
\cite{Blackwell56}. Note that it also follows from $(i).$ \qed

A straightforward application of this last theorem is given by the
following result. It will be used  several times in the forthcoming
sections.

Let $\mu \in \RR^N$ with $\mu \neq 0.$ For all $x \in \RR^N$, set
$\mu(x) = \langle \mu, x \rangle.$
\bcor \label{th:Bcor}
Suppose there exist actions $a_1, a_2 \in A$ and numbers $\alpha,
\beta$ such that for all $b \in B$
$$\mu (U(a_1,b)) \leq \alpha \mbox{ and }  \mu(U(a_2,b)) \geq \beta.$$
Let $Q$ be a payoff-based strategy such that
$$\mu(u) > \alpha \Rightarrow Q_u(a_1) = 1$$ and $$\mu(u)  < \beta
\Rightarrow Q_u(a_2) = 1.$$ Then
$$\Lambda = \{u \in E :  \mu(u) \in [\alpha, \beta] \}$$ is  a $B-$set
\ecor
Note that there is no assumption here that $\alpha \leq \beta.$ If
$\alpha \geq \beta$ $[\alpha, \beta]$ stands for $[\beta,\alpha].$

\prf Equation (\ref{eq:separ}) in this context becomes
 $$\mu(u) > \alpha  \Rightarrow \mu(v) \leq \alpha$$
 $$\mu(u) < \beta \Rightarrow \mu(v) \geq \beta,$$ for all $u \in E, v
\in C(u).$ By convexity of the half spaces $\{\mu(v) \leq \alpha \},
\{\mu(v) \geq \beta \},$ and the definition of $Q$ this is equivalent
to the condition given in the statement of the corollary.  \qed

We conclude this section with some quantitative estimates given in the
excellent recent survey paper by Perchet \cite{Perchet14}.
Let $$|E| = \sup \{ \|v\| \: : v \in E \}, \ |\Lambda| =  \sup \{ \|v\|
\: : v \in \Lambda \}.$$
 The first assertion of the next theorem
follows from Corollary 1.1 in \cite{Perchet14}. It is slight variant
of a result obtained by Blackwell \cite{Blackwell56}. The second
assertion follows from Corollary 1.5 in \cite{Perchet14}.
\bthm
\label{th:Black2}
Suppose DM adopts the payoff-based strategy $Q$ and that  $\Lambda
\subset E$ is a ${\cal B}$-set for $Q.$ Then for all $\eta > 0$

\bdes

\iti $$\Pr(\sup_{m \geq n} d(u_m, \Lambda) \geq \eta) \leq \frac{ 2
  (|E| + |\Lambda|)^2}{\eta^2 n}$$.

\itii If furthermore $\Lambda$ is convex,
$$\Pr(\sup_{m \geq n} d(u_m, \Lambda) - 2 \frac{|E|}{\sqrt{m}} \geq
\eta) \leq 4 \exp (- \frac{\eta^2 n}{32 |E|^2}).$$
\edes
\ethm

\section{$N$-Players Prisoner's Dilemma Game} \label{sec:NPD}
Consider an $N$-Players game (as described in Example
\ref{ex:Nplayers}) where  each player has two actions: {\em cooperate}
$C$ or {\em defect} $D$, so that $\Sigma^i = \{C,D\}.$ We assume that
the payoff functions $U^i, i = 1, \ldots N,$ are as follows. Let $s =
(s^1, \ldots, s^N) \in \{C,D\}^N$ be the action profile of the
players.
If Player $i$ cooperates (i.e $s^i = C$) and amongst her $N-1$
opponents, $k$ cooperate (i.e $\mathsf{card} \{j \neq i \: : s^j = C\}
= k$) she gets  $$U^i(s) =
v(C,k).$$ If she defects and amongst her  opponents, $k$ cooperate she
gets $$U^i(s) = v(D,k).$$ For $k = 0, \ldots N-1$ the numbers $v(C,k),
v(D,k)$ satisfy the following conditions, usual for prisoner's
dilemmas:

\bdes
\iti Defection is the dominant action :
$$v(C,k) < v(D,k)$$ for all $k = 0, \ldots, N-1;$
\itii The payoff of a defector increases with the number of cooperators :
$$v(D,k) \leq v(D,k+1)$$ for all $k = 0, \ldots, N-2;$
\itiii Mutual cooperation is a Pareto optimal:
For all $s \in \{C,D\}^N$
$$\sum_{i = 1}^N U^i(s) \leq N U^i(C, \ldots, C);$$ Or, equivalently,
for all $k = 0, \ldots, N-1$
$$k v(C,k-1) + (N-k) v(D,k) \leq N v(C,N-1).$$
\itiv (Occasional assumption) Mutual defection is Pareto inefficient:
For all $k = 0, \ldots, N-1$
$$k v(C,k-1) + (N-k) v(D,k) \geq N v(D,0).$$
\edes
\brem {\rm Condition $(i)$ makes $(D, \ldots, D)$ the unique Nash equilibrium of the one-shot game.} \erem
\brem {\rm
When $N = 2$ we retrieve the usual two players Prisoner's Dilemma
game. The Pareto conditions $(iii)$ and $(iv)$ amount to say that the
polygon with vertex set
 $$\{(v(D,0),v(D,0));\, (v(D,1),v(C,0));\, (v(C,0),v(D,1));\,
(v(C,1),v(C,1))\}$$
 is convex.}
\erem

\bex[Free riding]
{\rm Let $f : \{0, \ldots, N\} \mapsto \RR^+$ and $c > 0$ be such
  that $$\frac{c}{N} \leq f(k+1) - f(k) < c.$$
Let
$$v(C,k) = f(k+1) - c \mbox{ and } v(D,k) = f(k).$$  This can be seen
as a simple model of "free riding". Each player can either Contribute
(Cooperate), or Defect from contributing, to a public good. Individual
contribution costs $c$ and everyone -even if a defector- benefits from
the good and is paid $f(k)$, when there are $k$ contributors.

Note that the assumption on $f$ imply that conditions $(i)-(iv)$ above
are satisfied. The fact that mutual defection is a Nash equilibrium of
the one-shot game is reminiscent of  Hardin's book {\em The Tragedy of the
Commons}  \cite{HA68}.}
\eex

Recall that $E = \mathsf{conv}\{U(s) \: : s \in \{C,D\}^N\}.$
 For $u = (u_1, \ldots, u_N) \in E$ let
 $$\mu^i(u) = u_i - \frac{1}{N-1} \sum_{j \neq i} u_j = \frac{1}{N-1}
 \sum_{j = 1}^N (u_i - u_j).$$ Let $\delta$ be a nonnegative real
 number, Adapting \cite{Smale80}, \cite{BH96a} and \cite{BHS05}, we
 define a {\em $\delta$-good strategy} for Player $i$ as a
 payoff-based strategy $Q^i$ (as defined in section
 \ref{sec:notation}) for the Decision Maker, Player $i$, such that
 $$Q^i_u(C) = 1 \mbox{ if } \mu^i(u) \geq 0,$$ and
 $$Q^i_u(D) = 1 \mbox{ if } \mu^i(u) < - \delta.$$
 We call such a strategy {\em continuous} whenever the map $u \mapsto Q_u$ is continuous.

 The following result  shows that, by playing a $\delta$-good
 strategy, a player (or a group of players) makes sure that her
 opponents'average payoff cannot be much better than hers, nor than
 the Pareto optimal payoff. Under the supplementary condition $(iv)$
 she ensures that her payoff cannot be much worse that  the payoff
 resulting from mutual defection. If furthermore, all the players play
 a $\delta$-good strategy, one of them being continuous, the outcome
 is the one given by mutual cooperation.
  As a consequence (Corollary \ref{th:NPDNash}), continuous
  $\delta$-good strategies form a Nash equilibrium. The proof is
  postponed to the end of the section.

\bthm
\label{th:NPD} Let $k \leq N.$ Suppose that for all $i \in \{1,
\ldots, k\}$ Player $i$ plays a $\delta$-good strategy. Let
$$u_n^{-k} = \frac{1}{N-k} \sum_{j = k+1}^N u_n^j$$ be the average
payoff to players $k+1, \ldots, N$.  Then
\bdes
\iti
For all $i  \in \{1, \ldots, k\}$:
\[\begin{split}
&0 \leq \liminf_{n \rar \infty} u_n^{-k} - u_n^i \leq \limsup_{n \rar
  \infty} u_n^{-k} - u_n^i \leq \frac{N-1}{N-k} \delta,\\
&\limsup_{n \rar \infty} u^i_n \leq v(C,N-1),\\
\intertext{
and, if mutual defection is inefficient,}
& \liminf_{n \rar \infty} u^i_n  \geq v(D,0) - \delta.
\end{split}
\]
 \itii Suppose  $k = N$.  Then:
\bdes
\ita  $$L(\{u_n\}) \subset \mathsf{diag}(E) =
 \{u \in E \: : u_1 = \ldots = u_N\}, $$

and  if at least one of the players uses a continuous $\delta$-good strategy, then
\itb $$\lim_{n \rar \infty} u_n =  U(C, \ldots, C) = (v(C,N-1), \ldots v(C,N-1)).$$
\edes
\edes
\ethm

Let $\eps > 0.$  Let $\Theta^i$ be a strategy  (as defined by equation (\ref{eq:defstrat})) for Player $i.$

The strategy profile $(\Theta^1, \ldots,
\Theta^N)$
 is called an
      {\em $\eps$-Nash equilibrium} if for every $i$ and every alternative strategy  $\Xi^i$ for $i,$ the payoff to $i$ resulting from $(\Theta^1, \ldots, \Theta^{i-1}, \Xi^i, \Theta^{i+1}, \ldots, \Theta^N)$
        cannot be $\eps$ better than the payoff resulting from $(\Theta^1, \ldots,
\Theta^N).$ More precisely:

For every $i \in \{1, \ldots, N\},$ every strategy $\Xi^i$ and every
$\Sigma \times \Sigma$-valued process $\{(s_n, \tilde{s}_n)\}$ adapted
to the filtration $\{\F_n\}$ satisfying
$$\Pr(s_{n+1} = s | \cF_n) = \prod_{j = 1}^N \Theta_n^j(s^j),$$ and
$$\Pr(s_{n+1} = s | \cF_n) = \left ( \prod_{j \neq i} \Theta_n^i(s^j)\right )
\Xi_n^i(s^i);$$
then
$$\Pr(\limsup_{n \rar \infty} \tilde{u}_n^i \leq \liminf_{n \rar
  \infty} u_n^i + \eps) = 1.$$
Here $$u_n = \frac{1}{n} \sum_{k = 1}^n U(s_k), \tilde{u}_n = \frac{1}{n} \sum_{k = 1}^n U(\tilde{s}_k).$$
In other words, if all players but $i$ play the equilibrium strategy, Player $i$ cannot improve his payoff by more than $\eps$ if he deviates from $\Theta^i.$
\bcor
\label{th:NPDNash}
Let $\delta > 0.$ Suppose that for all $i \in \{1, \ldots, N\}$  $Q^i$
is a continuous $\delta$-good strategy.
Then $(Q^1, \ldots, Q^N)$ is a $\delta(N-1)$-Nash equilibrium.
\ecor
\prf Follows from Theorem \ref{th:NPD} \qed
\subsection*{Proof of Theorem \ref{th:NPD}}
For all $\delta \geq 0$ let
  $$\Lambda^i(\delta) = \{u \in E \: : - \delta \leq \mu^i(u) \leq 0\}.$$
\bprop
 \label{th:NPDapproach} Let $i \in \{1, \ldots, N\}.$
 Suppose that $i$ plays a $\delta$-good strategy $Q^i$. Then $\Lambda^i(\delta)$
 is a ${\cal B}$-set for $Q^i.$
 In particular, assertion $(ii)$ of Theorem \ref{th:Black} and Theorem  \ref{th:Black2} hold.
 \eprop
\prf Suppose $i = 1.$
Given $s^{-1} = (s^2, \ldots, s^N) \in \{C,D\}^{N-1},$ let $k = \mathsf{card}\{j > 1 \: : s^j = C\}.$
Then $$\mu^1(U(C,s^{-1})) = \frac{N-(k+1)}{N-1}(v(C,k) - v(D,k+1)) \leq 0$$
and
$$\mu^1(U(D,s^{-1})) = \frac{k}{N-1}(v(D,k) - v(C,k-1)) \geq 0.$$
By Corollary \ref{th:Bcor} and definition of $Q^1$, this concludes the proof. \qed
 \bprop[Properties of $\{\Lambda^i(\delta)\}$]
 \label{th:propLambda}
 For all $\delta \geq 0$ the sets $\Lambda^i(\delta)$ satisfy the
 following properties:
 \bdes
 \iti
 For all $i \in \{1, \ldots, N\}$ and  $u \in \Lambda^i(\delta)$
 $$u_i \leq v(C,N-1),$$ and, if mutual defection is Pareto inefficient,  $$u_i \geq v(D,0) - \delta.$$
 \itii For all $k \leq N,  u \in \bigcap_{i  \in \{1, \ldots, k\}} \Lambda^i(\delta)$ and $i \in \{1, \ldots, k\}$
$$0 \leq \frac{\sum_{j = k+1}^N u_j}{N-k}  - u_i  \leq \delta \frac{N-1}{N-k}$$
\itiii  $$ \bigcap_{i  \in \{1, \ldots, N\}} \Lambda^i(\delta) = \mathsf{diag}(E) =
\{u \in E \: : u_1 = \ldots = u_N\}.$$
 \edes
\eprop
\prf
Suppose $i = 1.$
For all $\eta \in \RR$ let $\Pi^{\eta}$ be the orthogonal projection onto the hyperplan $\{\mu^1(u) = \eta\}.$
Then $$\Pi^{\eta}(u) = u - (\frac{\mu^1(u) - \eta}{\|\mu^1\|^2})\mu^1$$
where $\mu^1$ is the vector defined by $\mu^1(u) = \langle \mu^1 , u \rangle.$ That is $\mu^1 = (\mu^1_i)_{i = 1, \ldots, N}$ with
$\mu^1_1 = 1$ and $\mu^1_i = -  \frac{1}{N-1}$ for $i > 1.$ It follows that
$$\Pi^{\eta}_1(u) = u_1 - \frac{\mu^1(u) - \eta}{\|\mu^1\|^2} = \frac{\sum_{i = 1}^N u_i}{N} + \eta \frac{N-1}{N}.$$
Thus, by Pareto dominance of $v(C,N-1)$ $$\Pi^{\eta}_1(u) \leq v(C,N-1) + \eta$$  for all $u = U(s),$ hence for all $u \in E.$
Let now $u \in \Lambda^1.$ Then $u = \Pi^{\eta}(u)$ with $\eta = \mu^1(u) \in [-\delta,0].$ Thus $u_1 \leq v(C,N-1).$
Similarly, if  $v(D,0)$ is inefficient, then $u_1 \geq v(D,0) - \delta \frac{N-1}{N}.$

To prove the second  assertion set $A = \sum_{j = 1}^k u_j, B = \sum_{j > k} u_j$ and note that, by definition of $\Lambda^i(\delta)$
$$-(N-1) \delta \leq N u_i - A - B \leq 0$$ for all $i = 1, \ldots, k.$ Thus, by summing over all $i = 1, \ldots, k,$
$$ -(N-1) k \delta \leq (N-k) A - k B \leq 0.$$ Then
$$ N u_i \leq A + B \leq  \frac{k B}{N-k} + B = \frac{N B}{N-k}$$
 and
 $$  N u_i \geq A + B -(N-1)\delta \geq \frac{k (B - \delta (N-1))}{N-k} + (B- (N-1) \delta)= \frac{N (B- (N-1) \delta)}{N-k} .$$
 The last assertion  is immediate, because on $\Lambda^i(\delta)$ $N u_i \leq \sum_j u_j.$
\qed
\paragraph{Proof of Theorem \ref{th:NPD}}
 Assertion $(i)$ and the beginning of $(ii)$ follow from Propositions \ref{th:NPDapproach} and \ref{th:propLambda}.
It remains to prove the last assertion. Assume that Player $1$ uses a continuous strategy. Recall that $C^1(x)$ is defined by (\ref{eq:defC}) with $Q = Q^1.$ By continuity of $x \mapsto Q^1_x$ the map $x \mapsto C^1(x)$ has a closed graph so that $\overline{\mathsf{co}}(C)^1(x) = C^1(x)$ for all $x \in E.$ Thus,
by Theorem \ref{th:limitset} $(i)$,   the set $L = L(\{u_n\})$ is invariant under the differential inclusion
\beq
\label{eq:difC1}
\dot{u} \in - u + C^1(u)
\eeq and, by what precedes, is contained in $\mathsf{diag}(E).$ In particular,
for all $u  \in L$ there exists $\eta$ solution to (\ref{eq:difC1})
such that $\eta(0) = u$ and  $\mu^1(\eta(t)) = 0$ for all $t.$ Let
$h(t) = \eta(t) + \dot{\eta}(t).$ Then $\mu^1(h(t)) = 0$ and $h(t) \in
C^1(\eta(t))$ for almost all $t.$
Let $v^* = U(C, \ldots, C) = (v(C,N-1), \ldots, v(C,N-1)).$

By definition of $Q^1$ and $C^1,$ $$\mu^1(u) \geq 0 \Rightarrow C^1(u)
= \mathsf{conv} \{U(C,s^{-1}) : s^{-1} \in \{C,D\}^{N-1}\}.$$ Now, the
proof of Theorem  \ref{th:NPDapproach} shows that $\mu^1(U(C,s^{-1})
\leq 0$   with equality only if $s = (C, \ldots, C).$ Thus  $$\mu^1(u)
\geq 0 \Rightarrow \{v \in C^1(u) \: : \mu^1(v) = 0\} = \{v^*\}.$$
This implies that $h(t) = v^*$ and $\eta(t) = e^{-t}(u- v^*) + v^*$
for all $t \in \RR.$ By compactness of $L$ we must have $u = v^*$ (for
otherwise $\{\eta(t)\}$ would be unbounded).

\section{Network Prisoner's Dilemma Games} \label{sec:graphPD}
\subsection*{Network Games}

In this section we consider a game in which players are located at the
vertices of a graph and interact only with their neighbors. There are
$M$ players denoted $i = 1, \ldots, M.$ Player $i$ has a finite action
set $\Sigma^i.$ The set $V = \{1, \ldots, M\}$ of vertices of the
graph is equipped with an {\em edge set} ${\cal E} \subset V \times
V.$

We assume that the graph $(V,{\cal E})$ is
 \bdes \ita symmetric:
$(i,j) \in {\cal E} \Rightarrow (j,i) \in {\cal E},$

\itb self-loop free: $(i,i) \not \in {\cal E},$ and

\itc  irreducible: for all $i,j \in V $
there exist $k \geq 1$ and $i_1, \ldots, i_k \in V$ such that $i_1 =
1, i_k = j$ and $(i_{l}, i_{l+1}) \in {\cal E}$ for $l = 1, \ldots,
k-1.$ \edes For each $(i,j) \in {\cal E}$ there is a real valued
map $$U^{ij} : \Sigma^i \times \Sigma^j \mapsto \RR$$ representing the
{\em payoff function to Player $i$ against Player $j.$}

Let $\mathsf{Neigh}(i) = \{j  \in V \: : (i,j) \in {\cal E} \}$ and let $N_i$ be its cardinal.
The {\em payoff function to $i$} is the map $U^i : \Sigma  \mapsto \RR^{N_i}$ defined by
$$U^i(s) = (U^{ij}(s^i,s^j))_{j \in \mathsf{Neigh}(i)}.$$
Using the notation of Example \ref{ex:Nplayers}, set $N = \sum_{i = 1}^M N_i,$ and define the vector payoff function of the game as
 $$U = (U^1, \ldots, U^M) : \Sigma \mapsto  \RR^{N_1} \times \ldots \times \RR^{N_m} \simeq \RR^N.$$ The state space of
the game is then $E = \mathsf{conv}\{U(s), s \in \Sigma\} \subset
\RR^N.$

In addition to these data, we assume given a Markov transition matrix $K = (K_{ij})_{i,j \in V}$ adapted to $(V, {\cal E}).$
That is
$$K_{ij} \geq 0, \quad \sum_j K_{ij} = 1$$
and
$$K_{ij} > 0 \Leftrightarrow (i,j) \in {\cal E}.$$
The {\em mean payoff} to Player $i$ for the strategy profile $s$ is defined as
\beq
\label{eq:defglob}
\overline{U}^i(s) = \sum_j K_{ij} U^{ij}(s).  \eeq Irreducibility of
the graph $(V,{\cal E})$ ensures irreducibility of the transition
matrix $K$.  Therefore there is a unique {\em invariant probability}
$\pi$ for $K.$ That is, $$\pi_i \geq 0, \sum_i \pi_i = 1$$ and for all
$i \in V$ $$\sum_j \pi_j K_{ji} = \pi_i.$$ Define the {\em weight of
  edge} $(i,j) \in {\cal E}$ as \beq
\label{eq:defomega}
\omega_{ij} = \pi_i K_{ij}.
\eeq
Such weights will prove to be useful for defining $\delta$-good strategies below.
Note that, by invariance of $\pi,$
\beq
\label{eq:omega}
\sum_j \omega_{ij} = \sum_j \omega_{ji} = \pi_i
\eeq
\bex
\label{ex:graph}
{\rm Suppose  $$K_{ij} = \left \{ \begin{array}{c}
                               \frac{1}{N_i}  \mbox{ if } j \in \mathsf{Neigh}(i)\\
                              0 \mbox{ if } j \not \in \mathsf{Neigh}(i)
                            \end{array} \right.$$
Then $$\overline{U}^i(s) = \frac{\sum_{j \in \mathsf{Neigh}(i)} U^{ij}(s)}{N_i},$$
$$\pi_i = \frac{N_i}{N} \mbox{ and } \omega_{ij} = \frac{1}{N} \Ind_{j \in \mathsf{Neigh}(i)}.$$
}
\eex
\subsection*{Network Prisoner's Dilemma Games}
We consider now a particular example of network games where each pair of neighboring players is engaged in two players prisoner dilemma game.
We assume that for each $i \in V$ $\Sigma^i = \{C,D\},$ and
$$U^{ij}(C,D)  = CD, \ldots, U^{ij}(D,C) = DC,$$
where
\bdes
\iti
$$CD < DD < CC < DC$$
as usual for the two player prisoner's dilemma game.
\itii We furthermore assume that the  outcome $CC$ is {\em Pareto optimal} and that the outcome $DD$ is {\em Pareto inneficient}, in the sense that for all $(i,j) \in {\cal E}$
$$(\omega_{ij} + \omega_{ji}) DD < \omega_{ij} CD + \omega_{ji} DC < (\omega_{ij} + \omega_{ji}) CC;$$
\edes
\brem {\rm If $K$ is reversible with respect to $\pi$ (meaning that
  $\omega_{ij} = \omega_{ji})$ as in Example \ref{ex:graph}, Pareto
  inefficiency means
$$2 DD < CD + DC < 2 CC.$$
Equivalently, the polygon with vertices
$$(DD,DD),\, (CD,DC),\,(DC,CD),\, (CC,CC)$$
 is convex and hence equal to $E$.}
\erem
For $u = (u_{ij})_{i \in V, j \in \mathsf{Neigh}(i)} \in \RR^{N_1} \times \ldots \times \RR^{N_m} $


let
$$\mu^i(u) = \sum_j \omega_{ij} u_{ij} - \omega_{ji} u_{ji}.$$
Given $\delta \geq 0,$ a {\em $\delta$-good strategy} for Player $i$ is a payoff-based strategy $Q^i$ such that
$$Q_u^i(C) = 1 \mbox{ if } \mu^i(u) \geq 0,$$ and
$$Q_u^i(D) = 1 \mbox{ if } \mu^i(u) < -\delta.$$

The following result is similar to Theorem \ref{th:NPD}. It shows that
if a group of players use $\delta$-good strategies, their payoffs
cannot be much worse that the payoff resulting from mutual defection
and that a weighted average of the other players payoffs cannot be
much better than hers.
If furthermore, all the players play a
$\delta$-good strategy, and that of player $i$ is
continuous, then the payoffs of $i$ against $j$  and $j$ against $i$
both equal
$CC$, given by mutual cooperation.

As a consequence (Corollary
\ref{th:NetPDNash}), continuous $\delta$-good strategies form a Nash
equilibrium. The proof is postponed to the end of the section.
 \bthm
\label{th:NetPD} Assume $1\le k \leq N.$
Suppose that for  $i \in \{1,\ldots, k\}$, Player $i$ plays a
$\delta$-good strategy. Then
\bdes

\iti
$$L(\{u_n\}) \subset \bigcap_{i  \in \{1, \ldots, k\}}
\Lambda^i(\delta).$$

\itii
$$DD - \frac{\delta}{2 \pi_i} \leq \liminf_{n \rar \infty}
\overline{u}_n^i \leq \limsup_{n \rar \infty} \overline{u}_n^i \leq
CC, \qquad (i=1,\dots,k).$$

\itiii
$$\sum_{j = k+1}^N \pi_j DD \leq \liminf_{n \rar \infty} \sum_{j = k+1}^N \pi_j \overline{u}^j_n \leq \limsup_{n \rar \infty} \sum_{j = k+1}^N \pi_j \overline{u}^j_n \leq  \sum_{j = k+1}^N \pi_j CC + \frac{k \delta}{2}.$$

\itiv If $k = N$ and Player $l$ uses a continuous $\delta$-good
strategy, then for all $j \in \mathsf{Neigh}(l)$
$$\lim_{n \rar \infty} u^{lj}_n = \lim_{n \rar \infty} u^{jl}_n= CC$$
\edes
\ethm

\bcor
\label{th:NetPDNash} Let $\delta > 0.$ Suppose that for all $i \in \{1, \ldots, N\}$ $Q^i$ is a continuous $\delta$-good strategy. Then $(Q^1, \ldots, Q^N)$ is a $\frac{(N-1) \delta}{2}$ Nash equilibrium.
\ecor
\subsection*{Proof of Theorem \ref{th:NetPD}}
For all $\delta \geq 0$ let $$\Lambda^i(\delta) = \{u \in E \: :
-\delta \leq \mu^i(u) \leq 0\}.$$
\bprop
\label{th:NetPDapproach} Let $i \in \{1,\ldots, N\}$.
Assume that $i$ plays a $\delta$-good strategy $Q^i.$ Then $\Lambda^i(\delta)$
 is a ${\cal B}$-set for $Q^i.$
\eprop
\prf Fix $i \in V.$
Let $s = (s^1, \ldots, s^M) \in \Sigma^1 \times \ldots \times \Sigma^M$ be such that $s^i = C.$ Then
$$\mu^i(U(s)) = \sum_j \omega_{ij} (CC t_j + CD (1-t_j)) - \sum_j \omega_{ji} (CC t_j + DC (1-t_j)) )$$
where $t_j = 1$ if $s^j = C$ and $0$ otherwise.
Thus
$$\mu^i(U(s)) = \sum_j \omega_{ij} t_j (CC - CD) + \sum_j \omega_{ji} t_j (DC - CC) + CD \sum_j \omega_{ij} - DC \sum_j \omega_{ji}$$
$$\leq \sum_j \omega_{ij} (CC -CD) + \sum_j \omega_{ji} (DC -CC) + CD \sum_j \omega_{ij} - DC \sum_j \omega_{ji}$$
$$ = \pi_i (CC - CD + DC - CC + CD - DC) = 0.$$
Suppose now that $s^i = D.$ Then
$$\mu^i(U(s)) =  \sum_j \omega_{ij} (DC t_j + DD (1-t_j)) - \sum_j \omega_{ji} (CD t_j + DD (1-t_j)) )$$
$$ =  \sum_j \omega_{ij} t_j (DC - DD) + \sum_j \omega_{ji} t_j (DD - CD) \geq 0.$$
The results then follows from Corollary \ref{th:Bcor}. \qed
\brem
\label{rem:NetPDapproach}
{\rm the proof above shows that $\mu^i(U(s)) < 0$ (respectively $> 0$) if $s^i = C$ (resp. $D$) and $s^j = D$ (resp. $C$) for some $j \neq i$ }
\erem
\bprop[Properties of $\{\Lambda^i(\delta)\}$]
\label{th:proplambda2}
For all $\delta \geq 0$ the sets $\{\Lambda^i(\delta)\}$ verify the
following properties:
\bdes \iti For $i \in \{1, \ldots, N\}$ and $u
\in \Lambda^i(\delta)$ set $\overline{u}_i = \sum_j K_{ij} u_{ij}$.
Then
$$DD - \frac{\delta}{2 \pi_i}  \leq \overline{u}_i \leq CC.$$
\itii For all $k < N$  and $u \in \bigcap_{j = 1}^k \Lambda^j(\delta)$
$$0 \leq \sum_{i = k+1}^N \mu^i(u) \leq k \delta,$$
$$(\sum_{i = k+1}^N \pi_i)  DD \leq \sum_{i = k+1}^N \pi_i \overline{u}_i \leq (\sum_{i = k+1}^N \pi_i) CC  + \frac{k \delta}{2} $$
\itiii $$\bigcap_{i \in \{1, \ldots, N\}} \Lambda^i(\delta) = \bigcap_{i \in \{1, \ldots, N\}} \Lambda^i(0)$$
\edes
\eprop
 \prf
$(i).$ Let $v \in \RR^N$ be the vector defined  by $v_{ij} = 1, v_{ji} = -1$ for all $j \in \mathsf{Neigh}(i)$
and $v_{kl} = 0$ if $k \neq i, l \neq i$ or $(k,l) \not  \in {\cal E}.$
Let $\Pi^{\eta}$ be the projection onto the hyperplan $\{\mu^i(u) = \eta\}$ parallel to $v.$ That is
$$\Pi^{\eta}(u) = u - \frac{\mu^i(u)-\eta}{\mu^i(v)} v = u - \frac{\mu^i(u)-\eta}{2 \pi_i} v.$$
Thus, for all $s \in \Sigma$
$$\sum_j \omega_{ij} \Pi^{\eta}(U(s))_{ij}  = \frac{\sum_j (\omega_{ij} U^{ij}(s) + \omega_{ji} U^{ji}(s)) + \eta}{2}$$
\begin{eqnarray}
\label{eq:pieta}
 & \leq& \frac{\sum_j (\omega_{ij}  + \omega_{ji}) CC  + \eta}{2} = \pi_i CC + \frac{\eta}{2}, \\
 \label{eq:pieta2}
 & \geq & \frac{\sum_j (\omega_{ij}  + \omega_{ji}) DD  + \eta}{2} = \pi_i DD + \frac{\eta}{2}
\end{eqnarray}
where the last inequalities follow from Pareto dominance. This implies
that for all $u \in E$
$$DD \frac{\eta}{2 \pi_i} \leq \sum_{j} K_{ij} \Pi^{\eta}(u)_{ij} =
\frac{1}{\pi_i} \sum_{j} \omega_{ij} \Pi^{\eta}(u)_{ij} \leq CC +
\frac{\eta}{2 \pi_i}.$$ Hence for all $u \in \Lambda^i$
$$DD - \frac{\delta}{2 \pi_i} \leq \sum_{j} K_{ij} u_{ij} \leq CC.$$

$(ii).$ Note that
$$\sum_j \mu^j(u) = \sum_{j } \sum_{k } \omega_{jk} u_{jk} -
\omega_{kj} u_{kj} = \sum_{k } \sum_{j } \omega_{jk} u_{jk} -
\omega_{kj} u_{kj} = - \sum_j \mu^j(u).$$ Thus, $\sum_{j} \mu^j(u) =
0$ and the inequalities follow from the definition of
$\Lambda^j(\delta)$ for the first one and inequalities
(\ref{eq:pieta}, \ref{eq:pieta2}) for the second one.

$(iii). $ Let $u \in \bigcap \Lambda^j(\delta).$ Then $\mu^i(u) \leq
0$ for all $i,$ but since $\sum_i \mu^i(u) = 0,$ $\mu^i(u) = 0.$
 \qed
 \paragraph{Proof of Theorem \ref{th:NetPD}} The proof is similar to
 the proof of Theorem \ref{th:NPD}.
 Assertion $(i)$, $(ii)$ and $(iii)$ follow from Propositions
 \ref{th:NetPDapproach} and \ref{th:proplambda2}. For $(iv)$ we use
 the fact that if player $1$ plays a continuous $\delta$ good
 strategy, then the limit set $L$ of $\{u_n\}$ is an invariant set of
 the differential inclusion $\dot{u} \in - u + C^1(u)$ contained in
 $\bigcap \Lambda^i(0).$ By proposition \ref{th:NetPDapproach} and
 remark \ref{rem:NetPDapproach}, for all $u \in \bigcap_{i}
 \Lambda^i(0)$ and $v \in C^1(u)$ $u^{1j} = u^{j1} = CC.$ Thus,
 reasoning like in the proof of Theorem \ref{th:NPD}, invariance of
 $L$ shows that for all $u \in L$ $u^{1j} = u^{j1} = CC.$
 \bibliographystyle{amsplain} \bibliography{prisoner}

\providecommand{\bysame}{\leavevmode\hbox to3em{\hrulefill}\thinspace}
\providecommand{\MR}{\relax\ifhmode\unskip\space\fi MR }
\providecommand{\MRhref}[2]{%
  \href{http://www.ams.org/mathscinet-getitem?mr=#1}{#2}
}
\providecommand{\href}[2]{#2}
\begin{thebibliography}{10}

\bibitem{Abhyankar01}
K.~Abhyankar, \emph{Smale strategies for prisoner's dilemma type games},
  Doctoral dissertation, University of California at Berkeley, 2001.

\bibitem{AX84}
R.~Axelrod, \emph{The evolution of cooperation}, Basic Book, Inc Publishers,
  New York, 1984.

\bibitem{B96}
M.~Bena\"{\i}m, \emph{A dynamical system approach to stochastic approximation},
  SIAM Journal on Optimization and Control \textbf{34} (1996), 437--472.

\bibitem{B99}
\bysame, \emph{Dynamics of stochastic approximation algorithms}, S\'{e}minaire
  de Probabilit\'es XXXIII, Lecture Notes in Math \textbf{1709} (1999), 1--68.

\bibitem{BH96}
M.~Bena\"{\i}m and M.~W. Hirsch, \emph{Asymptotic pseudotrajectories and chain
  recurrent flows, with applications}, J. Dynam. Differential Equations
  \textbf{8} (1996), 141--176.

\bibitem{BH96a}
M.~Bena\"{\i}m and M.W. Hirsch, \emph{Stochastic adaptive behavior for
  prisoner's dilemma}, Unpublished manuscript, University of California at
  Berkeley, 1996.

\bibitem{BHS05}
M.~Bena{\"\i}m, J.~Hofbauer, and S.~Sorin, \emph{{Stochastic approximations and
  differential inclusions}}, SIAM Journal on Optimization and Control
  \textbf{44} (2005), 328--348.

\bibitem{Ber78}
M.~Berger, \emph{G\'eom\'etrie, vol 3 : Convexes et polytopes, poly\`edres
  r\'eguliers, aires et volumes}, Fernand-Nathan, Paris, 1978.

\bibitem{Blackwell56}
D.~Blackwell, \emph{An analog of the minmax theorem for vector payoffs},
  Pacific Journal of Mathematics (1956), 1--8.

\bibitem{FR10}
M.~Faure and G~Roth, \emph{Stochastic approximations of set-valued dynamical
  systems: convergence with positive probability to an attractor}, Mathematics
  of Operation Research \textbf{35} (2010), 624--640.

\bibitem{HA68}
G.~Hardin, \emph{The tragedy of the commons}, Science (1968), 1234--1248.

\bibitem{Perchet14}
V.~Perchet, \emph{Approachability, regret and calibration: Implications and
  equivalences}, Journal of Dynamics and Games (2014), 181--253.

\bibitem{Smale80}
S~. Smale, \emph{The prisoner's dilemma and dynamical systems associated to
  non-cooperative games}, Econometrica (1980), 1617--1633.

\end{thebibliography}

\section*{Acknowledgments} We acknowledge support from the SNF grant  200020-149871/1.

\end{document}